\documentclass[11pt]{amsart}

\usepackage{amssymb,amsmath, amsfonts, amsthm, bbm, a4wide}

\newtheorem{theorem}{Theorem}[section]
\newtheorem{lemma}[theorem]{Lemma}
\newtheorem{corollary}[theorem]{Corollary}
\newtheorem{proposition}[theorem]{Proposition}

\newcommand{\G}{\varGamma}
\renewcommand{\L}{\varLambda}

\newcommand{\ts}{\hspace{0.5pt}}

\newcommand{\RR}{\mathbb{R}}
\newcommand{\QQ}{\mathbb{Q}}
\newcommand{\ZZ}{\mathbb{Z}}
\newcommand{\NN}{\mathbb{N}}
\newcommand{\HH}{\mathbb{H}}
\newcommand{\II}{\mathbb{I}}
\newcommand{\oo}{{\displaystyle\mathfrak{o}}}

\newcommand{\bt}{{\scriptscriptstyle \bullet}}

\newcommand{\ii}{\mathrm{i}}
\newcommand{\jj}{\mathrm{j}}
\newcommand{\kk}{\mathrm{k}}

\DeclareMathOperator{\lcm}{lcm}
\DeclareMathOperator{\den}{den}
\DeclareMathOperator{\tr}{tr}
\DeclareMathOperator{\nr}{nr}
\DeclareMathOperator{\N}{N}

\DeclareMathOperator{\OC}{\mathrm{OC}}
\DeclareMathOperator{\SOC}{\mathrm{SOC}}
\DeclareMathOperator{\OS}{\mathrm{OS}}
\DeclareMathOperator{\SOS}{\mathrm{SOS}}

\begin{document}

\title[Similar Sublattices and Coincidence Rotations of $A_4$ and its Dual]{Similar Sublattices and Coincidence Rotations \\[2mm] of the Root Lattice $\mathbf{A_4}$ and its Dual}

\author{Manuela Heuer} 

\address{Department of Mathematics and Statistics, The Open University, \newline
\hspace*{12pt}Walton Hall, Milton Keynes MK7 6AA, UK}
\email{m.heuer@open.ac.uk}

\begin{abstract}
 A natural way to describe the Penrose tiling employs the projection method on the basis of the root lattice $A_4$ or its dual. Properties of these lattices are thus related to properties of the Penrose tiling. Moreover, the root lattice $A_4$ appears in various other contexts such as sphere packings, efficient coding schemes and lattice quantizers.

Here, the lattice $A_4$ is considered within the icosian ring, whose rich arithmetic structure leads to parametrisations of the similar sublattices and the coincidence rotations of 
$A_4$ and its dual lattice. These parametrisations, both in terms of a single icosian, imply an index formula for the corresponding sublattices. The results are encapsulated in Dirichlet series generating functions. For every index, they provide the number of distinct similar sublattices as well as the number of coincidence rotations of $A_4$ and its dual.
\end{abstract}

\maketitle

\section{Introduction and generalities}
A \emph{lattice} $\G$ in  $\RR^d$, denoted by $\G=\langle b_1, \ldots, b_d \rangle_{\ZZ},$
 consists of all integer linear combinations of a basis $\{b_1, \ldots, b_d\}$ of $\RR^d$. In crystallography, coincidence site lattices of  a lattice $\G$ are used to describe and understand grain  boundaries; compare \cite{baake} and references given there. Mathematically, a \emph{coincidence site lattice} (CSL) of a lattice $\G$ is defined as 
$\G \cap R\G$, where $R$ is an isometry and the corresponding \emph{coincidence index} 

\[
   \varSigma (R) \, = \, [ \G : (\G\cap R\G)]\, 
\]
is finite, i.e.\! $R$ and $R\G$ share a common sublattice. Denote the orthogonal group by $\mathrm{O}(d,\RR)$ and define the set of all \emph{coincidence isometries} as
\[
    \OC (\G) \, := \,
    \{ R\in\mathrm{O}(d,\RR)\mid \varSigma (R) < \infty \}.
\] 
According to \cite{baake}, it is a subgroup of $\mathrm{O}(d,\RR)$. $\OC(\G)$ and $\SOC(\G)$, which is the subgroup of all coincidence rotations, can be characterised as subgroups of another important group of isometries, namely
\[
    \OS (\G) \, := \, 
    \{ R\in\mathrm{O}(d,\RR)\mid \alpha R\G\subset\G
    ,\alpha \in\RR, \alpha >0 \}.\, 
\] This group consists of all isometries that emerge form similarity mappings of the lattice $\G$ into itself; see \cite{csl} for details.  Its subgroup of rotations, denoted by $\SOS(\G)$, contains  $\SOC(\G)$ as a normal subgroup; see \cite{GB} for a detailed analysis of their relation. A sublattice of $\G$ of the form $\alpha R \G$ is called a \emph{similar sublattice} (SSL) of $\G$ and has obviously the index 
$
 [\G:\alpha R \G]=\alpha^d.
$
If the lattice $\G$ possesses a rich point symmetry, there are usually interesting SSLs, beyond the trivial ones of the form $m\G$ with $m \in \NN$.
For $R\in\OS (\G)$ define the \emph{denominator} relative to the lattice
$\G$ as
\begin{align*}  
   \den^{}_{\G} (R) \, := \,
   \min\ts \{\alpha \in  \RR \mid \alpha > 0,\; \alpha R\G\subset\G\}\, .
\end{align*}
As $R$ is an isometry, one always has $\den^{}_{\G} (R)\ge 1$,
and from $\den^{}_{\G} (R)\ts R\G\subset\G$ one concludes that
$\bigl(\den^{}_{\G} (R)\bigr)^d$ is an integer.  Consequently,
$\den^{}_{\G} (R)$ is either a positive integer or an irrational number,
but still algebraic and of degree $\leq d$.  Moreover, one has
\begin{align*} 
  \{ \alpha \in\RR\mid \alpha >0,\;\alpha R\G\subset\G\} \, = \,
  \den^{}_{\G} (R)\,\NN\, .
\end{align*}
This leads to the following characterisation of $\OC(\G)$ within $\OS(\G)$, compare \cite{csl} for details:
\begin{align} \label{OC-in-OS} 
\OC (\G) = \{ R\in \OS (\G)\mid \den^{}_{\G} (R) \in \NN \}. 
\end{align} 

When a lattice $\G \subset \RR^d$ is given, one is interested in the number of distinct SSLs of $\G$ of index $n$ as well as in the number of distinct CSLs of $\G$ of index $n$. If these arithmetic functions are multiplicative, they are often encapsulated into  Dirichlet series generating functions, because of their Euler product decomposition. A detailed introduction to arithmetic functions, the corresponding Dirichlet series and Euler products can be found in \cite{Apostol}. For many lattices in $d \leq 4$ the arithmetic functions which count the number of distinct SSLs and CSLs of each index have already been derived; see  for instance \cite{baake, BM98, BM99}. 
One lattice for which this problem has not been solved completely is the root lattice $A_4$. This lattice is of particular interest, because it forms the natural setting, in the sense of a minimal embedding,  for the description of the Penrose tiling as a cut and project set, see for example \cite{penrose}. Properties of the root lattice $A_4$ are thus directly related to the Penrose tiling. Some of the various other applications of the root lattice $A_4$ are described in \cite{CRS}.

In this paper recent results about the similar sublattices and coincidence rotations of the root lattice $A_4$ are reviewed and extended to its dual lattice $A^{*}_4$; details can be found in \cite{ssl, csl}.
In the following sections we first describe the root lattice $A_4$ and its dual $A_4^{*}$ in a suitable setting. Then we derive their Dirichlet series generating functions for the number of distinct SSLs of index $n$, which turn out to be the same for $A_4$ and $A_4^{*}$. For these lattices it is particularly difficult to determine the Dirichlet series generation functions for the number of distinct CSLs of index $n$. Therefore, we consider the slightly simpler problem to derive the generating functions for the number of distinct coincidence rotation of index $n$, which turn out to coincide, too.

\section{The root lattice $A_4$}
The root lattice $A_4$ is usually defined as  
$$A_4:=\{(x_1, \ldots, x_5) \in \ZZ^5 \mid x_1+\ldots+x_5=0\},$$ which lies in a $4-$dimensional hyperplane of $\RR^5$. One description of $A_4$ as a  lattice in $\RR^4$ is given by
\begin{align*}\label{def-alat}
L:=\bigl\langle (1,0,0,0), \tfrac{1}{2}&(-1,1,1,1), (0,-1,0,0), \tfrac{1}{2}(0,1,\tau\!-\!1,-\tau) \bigr\rangle_{\ZZ}\,
\end{align*}
where $\tau=(1+\sqrt{5})/2$ is the golden ratio. Relative to the inner product $\tr
(x\bar{y}) = 2 \langle x\ts | \ts y\rangle$,  where
$\langle x\ts | \ts y\rangle$ denotes the standard Euclidean
inner product, this lattice is the root lattice $A_4$; see \cite{CMP,ssl, csl} for
details. This particular description of the root lattice $A_4$ in $\RR^4$ is very convenient for our problem, as it enables us to use the arithmetic of the quaternion algebra $\HH (\QQ(\sqrt{5}\,))$; see \cite{KR} for a detailed introduction to Hamilton's quaternions. For brevity we use from now on the notation 
\[ 
K:=\QQ(\sqrt{5}\,)=\{q+r \sqrt{5} \mid  q,r \in \QQ \},
\] which is a quadratic number field.
The algebra $\HH (K)$ is explicitly given as
\[
\HH (K) = K \oplus \ii  K \oplus \jj  K \oplus \kk  K,
\]
 where the generating elements satisfy Hamilton's relations 
$\ii^2 = \jj^2 =
\kk^2 = \ii\jj\kk = -1.$ 
It is equipped
with a \emph{conjugation} $\,\bar{.}\,$ which is the unique mapping that
fixes the elements of the centre of the algebra $K$ and reverses the
sign on its complement. If we write 
\[
q=(a,b,c,d)=a+\ii b+\jj c+\kk d, \quad \text{this means}\quad \bar{q}=(a,-b,-c,-d).
\]
The reduced norm and trace in $\HH (K)$ are defined by
\[   \nr (q) \, := \, q\bar{q} \, = \, \lvert q\rvert^2 
   \quad \text{and} \quad
   \tr (q) \, := \, q + \bar{q}, \, \]
where we canonically identify an element $\alpha\in K$ with the
quaternion $(\alpha,0,0,0)$.  For any $q\in\HH (K)$, $\lvert q \rvert$ is
its Euclidean length, which need not be an element of $K$.
Nevertheless, one has $\lvert rs\rvert = \lvert r\rvert \lvert
s\rvert$ for arbitrary $r,s\in\HH (K)$.  Due to the geometric meaning,
we use the notations $\lvert q\rvert^2$ and $\nr (q)$ in parallel. 
An element $q\in\HH (K)$ is called \emph{integral} when both $\nr (q)$
and $\tr (q)$ are elements of \[\oo:=\ZZ[\tau]:=\{m+n\tau \mid m, n \in \ZZ \},\] which is the ring of
integers of the quadratic field $K$.

The \emph{icosian ring $\II$} consists of all linear combinations of the vectors 
\begin{align*} 
   (1,0,0,0),(0,1,0,0),
   \tfrac{1}{2}(1,1,1,1),\tfrac{1}{2}(1 \! - \!\tau,\tau,0,1)
   \end{align*} 
with coefficients in $\oo$. The elements of $\II$ are called \emph{icosians}. It is a remarkable object; see for example \cite{Reiner,MP,MW, ssl} and
references given there. In particular, $\II$ is a maximal order in $\HH (K)$, all elements of $\II$  are integral in $\HH (K)$ and the lattice $L$ is contained in $\II$.
Since $\HH (K)$ has class number $1$, compare \cite{Reiner, V}, all ideals
of $\II$ are principal. The detailed arithmetic structure of $\II$
is the key to the characterisation of the similar sublattices and coincidence rotations  \cite{ssl, csl} for
$L$. What is more, one significantly profits from another
map, called the \emph{twist map} in \cite{ssl}, which is an involution of
the second kind for $\HH (K)$; see \cite{KMRT} for details. If $q=(a,b,c,d)$, it is defined by the
mapping 
\[
q\mapsto \widetilde{q}:=\, (a{\ts}',b{\ts}',d{\ts}',c{\ts}')\, ,
\]
where ${}'$ denotes the algebraic conjugation in $K$, as defined by the
mapping $\sqrt{5}\mapsto -\sqrt{5}$. The algebraic conjugation in $K$ is
also needed to define the absolute norm on $K$, via 
$
     \N (\alpha) \, = \, \lvert \alpha \alpha{\ts}' \rvert\,  .
$
For the various properties of the twist map, we refer the reader to
\cite{ssl} and references therein.  The most important ones in our
present context are summarised in the following Lemma, whose proof can be found in \cite{ssl}. It describes the relations between $L$ and $\II$.
\begin{lemma} \label{fundamental}
Within $\HH (K)$, one has\/ $\widetilde{\II} = \II$ and
\[
L = \{ x\in\II\mid\widetilde{x}=x\}=\{ x + \widetilde{x}\mid x\in\II\} = \phi^{}_{+} (\II),
\] 
where the $\QQ$-linear
mapping $\phi^{}_{+} \! : \, \HH (K) \longrightarrow
\HH (K)$, is defined by $\phi^{}_{+} (x) = x+\widetilde{x}$.
\end{lemma}

The dual $A_4^{*}$ of the  the root lattice $A_4$, here in form of the dual of the lattice $L$, is given by
\[
L^{*}:=\{x \in \RR^{4} \mid \langle x|y \rangle \in \ZZ \text{ for all } y \in L\}.
\] 

\section{Similar Sublattices}
We are interested in the SSLs of the root lattice $A_4$ and its dual lattice $A_4^{*}$. According to \cite{ssl}, there is an index preserving bijection between the SSLs of $A_4^{*}$ and $A_4$, as well as between the SSLs of $A_4$ and $L$. Therefore, it is sufficient to concentrate on the SSLs of the lattice $L$.
For convenience, we define
$\HH (K)^\bt_{} = \HH (K)\setminus \{0\}$.
Proofs for this section can be found in \cite{ssl}.
\begin{lemma} \label{para-1}
   All SSLs of the lattice $L$
   are images of $L$ under orientation preserving mappings
   of the form $x\mapsto pxq$, with $p,q \in
   \HH (K)^\bt$.
\end{lemma}
\noindent This characterisation tells us that we only need to select an appropriate subset of \\ 
$\HH (K)^\bt \times \HH (K)^\bt$ in order to reach
all SSLs of $L$. A first step is provided by the following
observation.

\begin{lemma} \label{all-do}
    If $p\in\II$, $p L \widetilde{p}$ is an SSL of $L$.
\end{lemma}
On the other hand there is the following strengthening of Lemma \ref{para-1}.
\begin{proposition}  \label{fundamental-para}
   If $p L q \subset L$ with $p,q\in\HH (K)^\bt$,
   there is an $\alpha\in\QQ$ such that
 $
   q \; = \; \alpha\, \widetilde{p} \, .
$
\end{proposition}
To improve this characterisation we need the following definitions. An element $p\in\II$ is called \emph{$\II$-primitive} when
$\alpha p\in\II$, with $\alpha\in K$, is only possible with
$\alpha\in\oo$. Similarly, a sublattice $\L$ of $L$ is called \emph{$L$-primitive} when
$\alpha\L\subset L$, with $\alpha\in\QQ$, implies $\alpha\in\ZZ$. Whenever the context is clear, we simply use the term ``primitive'' in both cases.

\begin{corollary}  \label{para-2}
   All SSLs of the lattice $L$ are images of mappings of the
   form $x\mapsto \alpha \ts p x \widetilde{p}$ with $p \in\II$
   primitive  and $\alpha\in\QQ$. 
\end{corollary}

We now need to understand how an
SSL of $L$ of the form $p L \widetilde{p}$ with an
$\II$-primitive quaternion relates to the primitive sublattices of
$L$.

\begin{proposition}  \label{prim-one}
  If $p\in\II$ is\/ $\II$-primitive, $pL \widetilde{p}$ is an 
  $L$-primitive sublattice of $L$.
\end{proposition}

Combining the results of Corollary~\ref{para-2} and
Proposition~\ref{prim-one}, we obtain the following important
observation.
\begin{proposition} \label{prop-one-prime}
  A similar sublattice of $L$ is primitive if and only if
  it is of the form $pL\widetilde{p}$ with a primitive element $p \in \II$. Moreover, all SSLs of $L$ are of the form $q L \widetilde{q}$ with $q \in \II$. 
\end{proposition}

The next step is to find a suitable bijection that permits us to count the
primitive SSLs of $L$ of a given index.  Recall from \cite{MW,CMP}
that the unit group of $\II$ has the form
\begin{align*}
 \II^{\times} \, = \, \{ x\in\II \mid \N (\nr (x))= \pm 1 \},
 \end{align*} which leads to the following equivalence.

\begin{lemma}
  For $p\in\II$, one has $p L \widetilde{p} = L$
  if and only if $p\in\II^{\times}$.
\end{lemma}

Observe now that $p\ts\II=q\ts\II$ with $p,q\in\II$ holds if and only
if $q^{-1} p \in\II^{\times}$. The relevance of this fact in our
context comes from the observation that
\begin{equation*} 
    [\II : p\ts\II] \, = \, \N \bigl( \lvert p \rvert^4 \bigr)
    \, = \, [L : p L\widetilde{p}\,] \, ,
\end{equation*}
where the index of $p\ts\II$ in $\II$ follows from the determinant
formula in \cite[Fact 3]{ssl} followed by taking the norm in
$\ZZ[\tau]$.  This means that $p\ts\II \longleftrightarrow
p L\widetilde{p}$ describes an index preserving bijection between
primitive right ideals of $\II$ (meaning right ideals $q\II$ with
$q\in\II$ primitive) and primitive SSLs of $L$. 
Due to the definition of $L$-primitivity, a general SSL can be described as an integer multiple of a
primitive SSL. This leads together with the observation that all possible
indices of SSLs are squares, to the following central result.
\begin{theorem} \label{theoremssl}
    The number of SSLs of a given index $m$ is the same for the lattices
    $A_4^{*}, A_4$ and $L$. The possible indices are the squares of
    non-zero integers of the form $k^2 + k \ell - \ell^2$. If $f_{\rm SSL}(m)$ denotes the number of SSLs of index $m^2$, the corresponding Dirichlet series generating function reads
\[
 D_{\rm SSL} (s) \; := \; \sum_{m=1}^{\infty} \frac{f_{\rm SSL}(m)}{m^{2s}}
      \; = \; \zeta(4s)\, \frac{\zeta^{}_{\II}(s)}{\zeta^{}_{K} (4s)},
\]
where $\zeta(s)= {\displaystyle \prod_{p}}\frac{1}{1-p^{-s}}$ is the Riemann zeta function,
\begin{equation*}
\zeta^{}_{K} (s)= {\textstyle \frac{1}{1-5^{-s}}}\prod_{p \equiv \pm 2 (5)}{\textstyle \frac{1}{1-p^{-2s}}} \prod_{p \equiv \pm 1 (5)} {\textstyle \frac{1}{(1-p^{-s})^2}}
\end{equation*} 
 is the Dedekind zeta function of the quadratic field $K$, and $\zeta^{}_{\II} (s) \; = \; \zeta^{}_{K} (2s) \, \zeta^{}_{K} (2s-1)$
denotes the Dirichlet series for the right ideals of $\II$, which is the zeta function of $\II$. 
\end{theorem}  
Inserting the Euler products of $\zeta(s)$ and $\zeta^{}_{K} (s)$,
one finds the expansion of the Dirichlet series $D_{\rm SSL} (s)$ as an 
Euler product

\begin{align*}
 D_{\rm SSL} (s) \; = \; {\textstyle \frac{1}{(1-5^{-2s})(1-5^{1-2s})}}
  \hspace{-2mm} \prod_{p\equiv\pm 1\; (5)}\hspace{-1mm} {\textstyle \frac{1+p^{-2s}}{1-p^{-2s}}\,
   \frac{1}{(1-p^{1-2s})^2}} \prod_{p\equiv\pm 2\; (5)}
   {\textstyle \frac{1+p^{-4s}}{1-p^{-4s}}\,  \frac{1}{1-p^{2-4s}}} .
\end{align*}
Consequently, the arithmetic function $f_{\rm SSL}(m)$ is multiplicative and therefore completely specified by its values at prime powers $p^r$
with $r\ge 1$. These are given by 
\[
    f_{\rm SSL}(p^r) \; = \; \begin{cases}
        \frac{5^{r+1} - 1}{4} , & \text{$p=5$,} \\[1mm]
      \frac{2\, (1-p^{r+1}) - (r+1)(1-p^2)p^r}{(1-p)^2}, &
      \text{$p\equiv\pm 1 \; (5)$,} \\[1mm]
      \frac{2 - p^{r} - p^{r+2}}{1-p^2} , &
     \text{$p\equiv\pm 2 \; (5)$, $r$ even,} \\[1mm]
      0, &
     \text{$p\equiv\pm 2 \; (5)$, $r$ odd.}
    \end{cases}
\]
The first few terms of the Dirichlet series thus read

\[
   {\textstyle D_{\mathrm SSL} (s) \; = \; 1 + \tfrac{6}{4^{2s}} + \tfrac{6}{5^{2s}} + 
   \tfrac{11}{9^{2s}} + \tfrac{24}{11^{2s}} + \tfrac{26}{16^{2s}} + 
   \tfrac{40}{19^{2s}} + \tfrac{36}{20^{2s}} + \tfrac{31}{25^{2s}} +
   \tfrac{60}{29^{2s}} + \tfrac{64}{31^{2s}} + \tfrac{66}{36^{2s}}
   + \ldots}
\]
where the denominators are the squares of the integers previously
identified in \cite{CRS}.

\section{Coincidences and Rotations}
According to \cite[Th. 2.2]{baake}, $\OC(A_4)=\OC(A_4^{*})$ and the coincidence index of any conincidence isometry is the same for both lattices. Therefore, it is sufficient to investigate the CSLs of the lattice $A_4$. The best representation of this lattice is again the lattice $L$. First of all the investigation can be restricted to coincidence rotations without missing any CSLs, because $\overline{L}=L$, i.e. any orientation reversing operation can be obtained from an orientation preserving one after applying conjugation first.
According to Proposition~\ref{prop-one-prime}, a given SSL of $L$ is of the form $qL\widetilde{q}$, with $q \in \II$. The corresponding rotation is given by the mapping $x\mapsto
\frac{1}{\lvert q \widetilde{q}\ts\rvert}\,q x \widetilde{q}$.  Of course, many different icosians $q$ result in the same rotation.  Our aim is to restrict $q$ to suitable subsets of $\II$ without missing
any rotation.
In general, $\SOC(L)$ and $\SOS(L)$ are related according to Eq. (\ref{OC-in-OS}). Therefore, we have to identify the elements of $\SOC (L)$ within $\SOS(L)$. We would like to refer the reader to \cite{csl} for the proofs of this section.
\begin{lemma} \label{commensurate}
  Let $0\neq q\in\II$ be an arbitrary icosian.  The lattice
  $L\cap \frac{1}{\lvert q\widetilde{q}\ts\rvert}\, q L \widetilde{q}$
  is a CSL of $L$ if and only if\/ $\lvert q
  \widetilde{q}\ts\rvert\in\NN$.
\end{lemma}
Let us call an icosian $q\in\II$ \emph{admissible} when
$\lvert q \widetilde{q}\ts\rvert\in\NN$. As $\nr (\widetilde{q}\ts) =
\nr (q){\ts}'$, the admissibility of $q$ implies that
$\N\bigl(\nr (q)\bigr)$ is a square in $\NN$. With this definition the CSLs of $L$ can be characterised as follows.
\begin{theorem}
  The CSLs of $L$ are precisely the lattices of the form $L
    \cap \frac{1}{\lvert q\widetilde{q}\ts\rvert}\, q L\widetilde{q}$
    with $q\in\II$ primitive and admissible.
\end{theorem}
 This is the first step to define a bijection between certain primitive right
ideals $q\II$ of the icosian ring and the CSLs of $L$. The next step in this direction is the following Lemma.
\begin{lemma} \label{symmetries} 
  Let $r,s\in\II$ be primitive and admissible icosians, with
  $r\II=s\II$.  Then, one has 
\[L\cap \frac{r L\ts
    \tilde{r}}{\lvert r \tilde{r}\rvert} = L \cap \frac{s L
    \tilde{s}}{\lvert s \tilde{s}\rvert}.
\]
\end{lemma}
For our further discussion we need to replace the primitive and admissible icosian $p$, and with it $\tilde{p}$, by certain $\oo$-multiples, such that their norms have the
same prime divisors in $\oo$. In view of the form of the rotation $x\mapsto
\frac{1}{\lvert p\widetilde{p}\ts\rvert}\, p x \widetilde{p}$, this is actually rather natural because it
restores some kind of symmetry of the expressions in relation to
the two quaternions involved. For a primitive and admissible icosian $p\in\II$  we
choose explicitly 
\begin{align} \label{def-alpha}
    \alpha^{}_{q} \, = \, 
    \sqrt{\frac{\lcm(\nr (q),\nr (\widetilde{q}\,))}{\nr (q)}}
    \, \in \, \oo \, ,
\end{align}
where we assume a suitable standardisation for the $\lcm$ of
two elements of $\oo$; see \cite{csl} for details. Moreover, we have the relation $\alpha^{}_{\tilde{q}}
= \widetilde{\alpha^{}_{q}} = \alpha^{\,\prime}_{q}$.
The icosian $\alpha_{q} q$ is called the \emph{extension} of the
primitive admissible element $q\in\II$, and $(\alpha^{}_q q,
\alpha^{\,\prime}_{q} \widetilde{q}\,)$ the corresponding
\emph{extension pair}.  Since $\alpha_q$ and $\alpha_q'$ are central, the extension does not change the rotation, i.e.
\begin{align*} 
     \frac{q x \widetilde{q}}{\lvert q\widetilde{q}\ts\rvert} 
     \, = \, \frac{q_{\alpha} x \widetilde{q_{\alpha}}}
     {\lvert q_{\alpha}\widetilde{q_{\alpha}}\ts\rvert}
\end{align*}
holds for all quaternions $x$.  Note that the definition of the
extension pair is unique up to units of $\oo$, and that one has the
relation
\begin{align*} 
     \nr (q_{\alpha}) \, = \,
     \lcm \bigl(\nr (q), \nr(\widetilde{q}\,)\bigr)
     \, = \, \nr (\widetilde{q_{\alpha}})
     \, = \, \lvert q_{\alpha}\ts
             \widetilde{q_{\alpha}}\rvert
     \, \in \, \NN \, ,
\end{align*}
which is crucial in the proof of the following Theorem.
{}For the further characterisation of the CSLs of $L$, it is convenient to define the set
\begin{equation}\label{sublattices} 
     L (q) \, = \, \{ qx + \widetilde{x}\widetilde{q}
     \mid x\in\II \} \, = \, \phi^{}_{+} (q\ts\II)\, ,
\end{equation}
which is a sublattice of $L$, compare Lemma~\ref{fundamental}.  Note
that, due to $\widetilde{\II}=\II$, one has $L(q) =
\widetilde{L(q)}$.
\begin{theorem} \label{csl-form}
   Let $q\in\II$ be admissible and primitive,
   and let $q_{\alpha} = \alpha_{q}\ts q$ be its extension.
   Then, the CSL defined by $\ts q$ is given by
\[
   L \cap \frac{1}{\lvert q\widetilde{q}\ts\rvert}\ts
    q L\widetilde{q} \, = \, L (q_{\alpha})\, ,
\]
   with $L (q_{\alpha})$ defined as in Eq.~$\eqref{sublattices}$.
\end{theorem}
\noindent This explicit identification of the CSL provides access to the corresponding index

\begin{theorem}  \label{index-formula}
  If $q\in\II$ is an admissible primitive icosian, the rotation
  $x\mapsto \frac{1}{\lvert q\widetilde{q}\ts\rvert}
  qx\widetilde{q}$ is a coincidence isometry of $L$.
  Moreover, the corresponding coincidence index is
\[
   \varSigma(q) \, = \, \nr(q^{}_{\alpha}) \, = \,
   \lcm (\nr(q),\nr(q)') \, ,
\]
which is, with our above convention from Eq.~$\eqref{def-alpha}$,
always an element of $\ts\NN$.  
\end{theorem}
 At this point, it is possible to determine the number  of distinct coincidence rotations of index $n$, where $n$ is a prime power. The rotations come in multiples of $120$, the order of the rotation symmetry group of $A_4$. 

\begin{theorem} The number of coincidence rotations of a given index $n$ is the same for the lattices
    $A_4^{*}, A_4$ and $L$. Let $120\ts f_{\rm SOC}(n)$ be the number of coincidence rotations of
  index $n$. Then, $f_{\rm SOC}(n)$ is a multiplicative
  arithmetic function, given at prime powers $p^r$ with $r \geq 1$ by
\begin{align*} \label{mult-fun}
   f_{\rm SOC}(p^r) \, = \, \begin{cases}
     6\cdot 5^{2r-1} , &  p=5  \, , \\
     \frac{p+1}{p-1}\, p^{r-1} (p^{r+1} + p^{r-1} - 2)
      , &  p\equiv\pm 1\; (5)\, , \\
     p^{2r}+p^{2r-2} , &  p\equiv\pm 2\; (5)  \, .
   \end{cases}
\end{align*} Its Dirichlet series generating function reads
\begin{eqnarray*} 
\lefteqn{D_{\rm \, SOC} (s)  \, =  
\sum_{n=1}^{\infty}\frac{f_{\rm SOC}(n)}{n^s} \,
 = \, \frac{\zeta^{}_{K} (s-1)}{1+5^{-s}}\, \frac{\zeta(s)\,\zeta(s-2)}{\zeta(2s)\,\zeta(2s-2)}} \\
   & &\, = {\textstyle \, \frac{1+5^{1-s}}{1-5^{2-s}}\,} \prod_{p\equiv\pm 1 \,(5)}
    {\textstyle \frac{(1+p^{-s})\,(1+p^{1-s})}{(1-p^{1-s})\,(1-p^{2-s})}}\,
    \prod_{p\equiv\pm 2 \,(5)} {\textstyle \frac{1+p^{-s}}{1-p^{2-s}}}  \\
   & &\, = \,{\textstyle 1 + \frac{5}{2^s} + \frac{10}{3^s} + \frac{20}{4^s} + 
      \frac{30}{5^s} + \frac{50}{6^s} + \frac{50}{7^s} + \frac{80}{8^s}+ \frac{90}{9^s} +  \frac{150}{10^s} + \frac{144}{11^s} +  \ldots},
\end{eqnarray*}
where $\zeta(s)$ is the Riemann zeta function and $\zeta^{}_{K} (s)$
  denotes the Dedekind zeta function of the quadratic field 
  $K$; see Theorem \ref{theoremssl} for their explicit expressions.
\end{theorem}
\section{Outlook}
Due to the factorisations of icosians into irreducible elements, which are difficult to access, we have not yet found a way to approach the number of CSLs of $A_4$ of index $n$ systematically. Nevertheless, it is possible to derive the number of CSLs up to index $120$, see \cite{csl}. We hope to report on some progress soon.

\section*{Acknowledgements}
It is a pleasure to thank M. Baake, U. Grimm, R.V. Moody and P. Zeiner for their suggestions and helpful discussions. This work was supported by EPSRC grant EP/D058465/1. Moreover, I thank the IUCr for their financial support to attend the ICQ10.

\end{document}